\documentclass[12pt]{amsart}

 \usepackage{amsfonts,amssymb,eucal}
\usepackage{amsthm}
 \usepackage{amsmath}
\usepackage{amscd}
 \usepackage{latexsym} 
\numberwithin{equation}{section}
\newcommand{\ad}{{\mbox{\rm ad}}}
 \newcommand{\e}{{\mbox{\rm e}}} 
 \newcommand{\li}{\rule[0.4ex]{1cm}{0.5mm}}
 \newcommand{\di}{\hspace*{1cm}}
 \newcommand{\mb}[1]{{\mbox{\boldmath{$#1$}}}}% mathematical bold
  \newcommand{\mc}[1]{{\mathcal{#1}}}% mathematical caligraphic
 \newcommand{\got}[1]{{\mathfrak{#1}}}% gothic with mbox for  mathematic
\newcommand{\db}[1]{{\mathbb{#1}}}% double
\newcommand{\gata}{\square} 
\newcommand{\pa}{\partial}
\newcommand{\R}{\ensuremath{\mathbb{R}}}
\newcommand{\C}{\ensuremath{\mathbb{C}}}
\newcommand{\N}{\ensuremath{\mathbb{N}}}

\newcommand{\T}{\ensuremath{S(U(1)\times U(1) \times U(1))}}
 \newcommand{\Hi}{\ensuremath{\mathcal{H}}}% Hilbert space
\newcommand{\sHi}{\ensuremath{{\mathcal{H}}^{*}}}% dual Hilbert space

\newcommand{\Hinf}{\ensuremath{\mathcal{H}^{\infty}}}% \infty Hilbert space
\newcommand{\g}{\ensuremath{\got{g}}}% Lie algebra g
\newcommand{\m}{\ensuremath{\got{m}}}% Lie algebra m
\newcommand{\bb}{\ensuremath{\got{b}}}% Lie algebra b
\newcommand{\gc}{\ensuremath{\got{g}_{\C}}}% Lie algebra g-complexificat
\newcommand{\Ugc}{\ensuremath{\mathcal{U}({\gc}})}% Universal algebra of
\newcommand{\Ad}{\ensuremath{{\mbox{\rm{Ad}}}}}% Ad- representation
\renewcommand{\P}{\ensuremath{\mathbb{P}}}
 \newcommand{\Ph}{\ensuremath{\P (\Hi )}}% Projective of Hilbert
\newcommand{\Phinf}{\ensuremath{\P (\Hinf )}}% Projective of Hilbert  % infty
\newcommand{\DM}{\ensuremath{{\got{D}}_M }}% D_M- module
\newcommand{\AM}{\ensuremath{{\got{A}}_M }}% A_M- module
 \newcommand{\D}{\ensuremath{{\got{D}}}}% D-module% now called sheaf
\newcommand{\A}{\ensuremath{{\got{A}}}}% A-module
\newcommand{\AAA}{\ensuremath{{\db{A}}_M}}% A-module
\newcommand{\am}{\ensuremath{{\bf{A}}_M}}% 
 \newcommand{\wt}{\widetilde}
\newcommand{\FSB}{symmetric Fock space }
 \newtheorem{Theorem}{Theorem}
 \newtheorem{Remark}{Remark}

\newcommand{\fl}{\ensuremath{{\mathcal{F}}_{\Hi}}}%Bargmann H space
\newcommand{\NC}{\ensuremath{\mathcal{G}}}% g in metric

\newcommand{\tc}{\ensuremath{\tilde{{\mb{C}}}}}% C tilde
\newtheorem{Proposition}{Proposition}
 \newtheorem{lemma}{Lemma}
\theoremstyle{definition}%nou 

\newtheorem{deff}{Definition}%martin

\begin{document}

\title{Realization of coherent state Lie algebras by 
 differential operators }
\author{Stefan  Berceanu}
 \address[Stefan  Berceanu]{National 
 Institute for Physics and Nuclear Engineering\\
         Department of Theoretical Physics\\
         PO BOX MG-6, Bucharest-Magurele, Romania}
\email{Berceanu@theor1.theory.nipne.ro}
\begin{abstract}
 A realization of
coherent state Lie algebras  by first-order differential operators
with holomorphic polynomial coefficients on  K\"ahler coherent state
orbits is presented.
 Explicit
formulas involving the Bernoulli numbers and the structure constants
for the semisimple Lie groups are proved. 
\end{abstract}
\maketitle
\noindent

\tableofcontents %%% ???? You want this?????
\section{Introduction}

The starting point of this investigation
is  the standard  Segal-Bargmann-Fock \cite{bar}
 realization
${\mb{a}}\mapsto \frac{\partial}{\partial z};~~ {\mb{a}}^+\mapsto z$ of the
canonical commutation relations $[{\mb{a}},{\mb{a}}^+]=1$ on the
symmetric Fock space $\fl :=\Gamma^{\mbox{\rm{hol}}}(\C,\frac{i}{2
\pi}\exp(-|z|^2)dz\wedge d\bar{z})$  attached to the Hilbert space
$\Hi :=L^2(\R, dx)$.  
The Segal-Bargmann-Fock realization can be considered as
a representation by differential operators of the real 3-dimensional
 Heisenberg algebra 
  $\g_{HW}=
< is1 + z{\mb{a}}^+ -\bar{z}{\mb{a}} >_{s\in\R;z\in\C}$ of
 the Heisenberg-Weyl group HW. We can look at this construction from
 group-theoretic point of view, 
considering  the complex number $z$  
 as local coordinate on the homogeneous manifold 
 $M:=HW/\R\cong \C$.  Glauber \cite{gl} has attached
 field coherent states (CS) to the points of the manifold $M$.   

We shall consider instead Glauber's \cite{gl} field CS
 generalized CS in the sense of Perelomov  \cite{perG} based on homogeneous
manifolds $M=G/H$.
 We  restrict ourself  to K\"ahler homogeneous spaces $M=G/H$
associated  to the so called CS-groups $G$,
 see \cite{mosc,mv,lis1,lis2,lis} and several works of Neeb quoted in
\cite{neeb}. The CS-groups
are groups whose quotients with  stationary groups are manifolds
which admit a holomorphic embedding in a projective Hilbert
space. This class of groups contains all compact groups, all simple
hermitian groups, certain solvable groups and also mixed groups as the
semidirect product of the Heisenberg group and the symplectic group
\cite{neeb}. We are interested in the
  realization of the  CS-Lie algebras by first order holomorphic
differential operators with polynomial coefficients.

The present work extends our previous results \cite{sbcag,sbl}.
 The differential action of the
 generators of the groups  on  coherent state manifolds
  which have the structure of
  hermitian symmetric spaces can be written down as a sum
 of two terms, one a polynomial $P$, and the second one  a sum of partial
 derivatives times some  polynomials $Q$-s,
    the degree of  polynomials  being  less than
 3 \cite{sbcag,sbl}.
It is interesting to investigate the same problem as in
 \cite{sbcag,sbl} on flag manifolds 
  \cite{bfr}.  Some
  results are available \cite{dob}, but they are not easily  handled. 
 We  give
 explicit formulas of the polynomials $P$ and $Q$-s  in the case of
semisimple Lie groups  
and also the simplest example of the compact nonsymmetric space
 $SU(3)/\T$, where the degree of the
polynomials is already 3.

 The paper is laid out as follows. \S \ref{CSR} collects some more of
less known facts about CS-groups and CS-representations. We follow 
\cite{lis1,lis2,lis} 
and  \cite{neeb}. Many of the facts summarized in \S \ref{CSR}
 have been already detailed
 in \cite{sbctim} and \cite{last}. 
The definition of  CS-groups is contained in \S  \ref{CSR1}.
\S \ref{CSvectors} defines the so called Perlomov's generalized
coherent state vectors in the context of the CS-groups.
 In \S \ref{fock} we recall the construction of the symmetric Fock
 space of
functions  \fl~  on which the differential operators  act.
In \S \ref{DIFF} we construct the representations of Lie algebras of
CS-groups by differential operators. \S\ref{mult}
recall some known facts about multipliers in the context of coherent
states. Data on  hermitian representations and differential operators are
summarized in \S\ref{sheaf}. Simple examples are presented in
\S \ref{SE}: the Heisenberg-Weyl group and $\got{sl}(2,\C )$.
 \S \ref{sSIMPLE} is dedicated to
the semisimple case.  The construction of Perelomov's coherent state
vectors
in this case is contained in \S\ref{perv}.
Our main new results are contained in Theorem
\ref{bigtheorem} in \S \ref{MAIN}. The simplest compact and
non-compact
 nonsymmetric
examples are contained in \S \ref{EXA}: $SU(3)/\T$ and $Sp(3,\R
)/\T$.

In this paper we give  the proof of the formulas
referring to the semisimple case. These formulas contain the Bernoulli
numbers and the structure constants. Let us recall that the Bernoulli
numbers appear \cite{kat} also in  connection  with
 Kontsevich's universal formula for deformation
quantization \cite{kon},  
 in the context of the Duflo-Kirillov 
isomorphism \cite{du} and Kashiwara-Vergne conjecture \cite{kw}.
Applications of the formulas here proved 
to explicit boson expansions
 for collective models on K\"ahler CS-orbits 
have been given in \cite{sbcbuc}. Using the same formulas of the
differential action of the generators of semisimple Lie groups, we have
written down the equations of motion generated by linear generators of
the groups on CS-orbits \cite{sbctim}. We have not included in this text
the holomorphic representation of CS-groups of semidirect product type
presented at the Conference {\it Operator algebras and Mathematical
Physics} in Sinaia, see \cite{jacobi}.

We have underlined  the deep
 relationship between coherent states and geometry \cite{progr}. Here
we are interested in the algebraic aspect.  
Our approach is closely related with those of
 reference  \cite{anton}, where are considered differential operators
 acting on coherent states
 constructed on Lie algebras. 

Let us also mention the  ``reflection symmetry''   approach \cite{jo}
  to study simultaneously 
the representations of the HW group and semisimple Lie groups.
  Applying the ``restriction principle'' as a particular
case it is obtained the Segal-Bargmann-Hall transform, 
 a generalization of the standard 
Segal-Bargmann
transform
 for compact groups, see references in \cite{hall}.

We want to underline that in this paper we do not give an explicit
construction of the symmetric Fock space of functions \fl~ on which
the differential operators act. For hermitian symmetric spaces a case
by case investigation was started by Hua \cite{hua} and developed 
by many authors, see e.g. \cite{tak,berezin2,ver,wal,john,up,far,iuri}.

We use for the scalar product the convention:
$(\lambda x,y)=\bar{\lambda}(x,y)$, $ x, y\in\Hi ,\lambda\in\C $.

\section{CS-representations, CS-vectors,  
reproducing kernel  }\label{CSR}
Perelomov's paper \cite{per} generated the interest of theoretician physicists
(e.g. \cite{onofri}) and mathematicians (e.g. \cite{mosc,mv}) in
understanding the mathematical aspects of CS-representations. Here we
follow the formulation of Lisiecki \cite{lis1,lis2,lis} and Neeb 
\cite{neeb}. The whole section  just fixes the definition and collects
some more or less known facts about CS-representations. More details
and proofs can be found in \cite{sbctim,last}.

\subsection{Coherent state representations}\label{CSR1}

Let us consider the triplet $(G, \pi , \Hi )$, where $\pi$ is
 a continuous, unitary
representation 
 of the  Lie group $G$
 on the   separable  complex  Hilbert space \Hi .
Let us denote by $\Hinf$ the {\it smooth vectors}. 
 Let us
pick up $e_0\in \Hinf$ and let  the notation:
$e_{g,0}:=\pi(g).e_0, g\in G$.
We have an action $G\times \Hinf\rightarrow\Hinf$, $g.e_0 :=
e_{g,0}$. When there is no possibility of confusion, we write just
$e_{g}$ for $e_{g,0}$. 

Let us denote by  
$[~]:\Hi^{\times}:=\Hi\setminus\{0\}\rightarrow\Ph=\Hi^{\times}/ \sim$ 
the projection with respect to the equivalence relation
 $[\lambda x]\sim [x],~ \lambda\in \C^{\times},~x\in\Hi^{\times}$. So,
$[.]:\Hi^{\times}\rightarrow \Ph , ~[v]=\C v$. The action 
$ G\times \Hinf\rightarrow\Hinf$ extends to the action
$G\times \Phinf\rightarrow\Phinf ,~ g.[v]:=[g.v]$. 

Let us now denote by $H$  the isotropy group $H:=G_{[e_0]}:=
\{g\in G|g.e_0\in\C e_0\}$.
We shall consider (generalized) coherent 
 states on complex  homogeneous manifolds $M\cong
G/H$, imposing the restriction that $M$ be  a complex submanifold of
\Phinf .

%\newpage
\begin{deff}\label{def1}
 a) The orbit $M$ is called a CS-{\em orbit} if 
there exits  a holomorphic embedding
$\iota : M \hookrightarrow \Phinf$.
In such a case $M$ is also called CS-{\it manifold}.

b) $(\pi,\Hi )$ is called a  CS-{\em representation} if there exists a cyclic
 vector $0\neq
e_0 \in\Hinf $ such that $M$ is a CS-orbit.

c)  The groups $G$ which admit
CS-representations  are called CS-{\it groups}, and their Lie algebras
\g~ are called CS-{\it Lie algebras}.
\end{deff}

For $X\in\g$, where \g ~is the Lie algebra of the Lie group $G$,  let
us define the unbounded operator $d\pi(X)$ on \Hi~ by
$d\pi(X).v:=\left. {d}/{dt}\right|_{t=0} \pi(\exp tX).v$,
whenever the limit on the right hand side exists. {\it The operator
$d\pi (X)$ is the closure of its restriction to \Hinf. In particular,
$id\pi (X)|_{\Hinf}$ is essentially selfadjoint}
 (cf. Proposition X.1.5. p. 391 in \cite{neeb}).
   We obtain a
representation of the Lie algebra \g~ on \Hinf , {\it the derived
representation}, and we denote
${\mb{X}}.v:=d\pi(X).v$ for $X\in\g ,v\in \Hinf$. Extending $d\pi$ by complex
linearity, we get a representation of the universal enveloping algebra
of the complex Lie algebra \gc~ on
the complex vector space \Hinf  
\begin{equation}\label{derived}
d\pi:\Ugc\rightarrow B_0(\Hinf ),~~ \mbox{\rm{with}} ~~ d\pi(X).v
:=\left.\frac{d}{dt}\right|_{t=0}\pi(\exp tX).v, X\in\g .
\end{equation}

   $B_0(\Hi^0)\subset \mc{L}(\Hi )$, where $\Hi^0:=\Hinf$, denotes the set
  of linear operators
 $A:\Hi^0\rightarrow \Hi^0$
which have  a formal adjoint  $A^{\sharp}:\Hi^0\rightarrow\Hi^0$, i.e.   
$(x,Ay)=(A^{\sharp}x,y)$ for all $x,y\in \Hi^0$.
Note that if  $B_0(\Hi^0)$ is the set of unbounded operators on \Hi ,
then the domain $\mathcal{D}(A^*)$ contains $\Hi^0$ and $A^*\Hi^0\subseteq 
 \Hi^0$, and  it makes sense to refer to the closure 
 of $A\in B_0(\Hi^0)$ (cf. \cite{neeb} p. 29; here $A^*$ is the
 adjoint of $A$).

We denote by $B:=<\exp_{G_{\C}} \bb >$ the Lie group
 corresponding to the Lie algebra
\bb ,  with
$\got{b}:=\overline{\got{b}(e_0)}$, where $\got{b}
 (v):=\{X\in\gc :X.v\in \C v\}=(\g_{\C})_{[v]}$. The group $B$ is
 closed in the complexification $G_{\C}$ of
$G$, cf. Lemma XII.1.2. p. 495 in \cite{neeb}. 
 The  complex structure on $M$ is induced by an 
embedding in a complex manifold,
 $i_1:M\cong G/H \hookrightarrow  G_{\C}/B$.
 We  consider such
manifolds  which admit a holomorphic
 embedding $i_2: G_{\C}/B\hookrightarrow
\Phinf$.   Then the embedding $\iota
=i_2\circ i_1$,  $\iota : M  \hookrightarrow \Phinf$
 is a holomorphic embedding,  and  the complex structure comes as in 
  Theorem XV.1.1 and Proposition XV.1.2 p. 646 in \cite{neeb}.
% So, the  $G$-invariant complex structures on the homogeneous
% space $M=G/H$   
% is introduced  in an algebraic manner. 

We conclude this paragraph recalling some known facts about CS-orbits
and CS-representations (cf. \cite{neeb}). 
Firstly, note that: {\it If $G[v]$ is a CS-orbit, then $v$ is an analytic
vector} (cf. Prop. XV 2.2 p. 651 in \cite{neeb}). 
Now, let  $G$ be a connected Lie group such that $\got{g}$ contains a
compactly embedded Cartan algebra $\got{t}$. 
Choosing $e_0:=v_{\lambda}$, where $\lambda \in i\got{t}^*$
 is a primitive element of a unitary
highest weight representation $(\pi_{\lambda},\Hi_{\lambda})$, then:
 $G.[v_{\lambda}]$ {\it is a complex orbit in bijection via the
momentum map  with the coadjoint
orbit} ${\mc{O}}_{-i\lambda}\subset\g^*$, and  
{\it every unitary highest weight representation
is a  CS-representa\-tion} (cf. Proposition XV.2.6 p. 652 in \cite{neeb}).
  {\em  The CS-representations of connected Lie groups are
irreducible} (cf. Proposition XV.2.7 p. 652 in
\cite{neeb}). Conversely,
{\it if \g~ is an admissible Lie algebra, $G$ a connected Lie group
with Lie algebra \g , and $(\pi,\Hi )$ a CS-representation of $G$ with
discrete kernel, then $\pi$ is a unitary highest weight
representation. If $G.[e_0]$ is a CS-orbit, then $e_0$ is a primitive
element for an appropriate positive system of roots with respect to a
compactly embedded Cartan subalgebra  and the orbit is the unique
complex orbit in \Phinf} (cf. Theorems XV.2.10 - 11 p. 655 in
\cite{neeb}).

\subsection{Coherent state vectors}\label{CSvectors}

Now we construct what we  call  Perelomov's
generalized  coherent state vectors, or simply CS-vectors,  based on the
CS-homogeneous manifolds $M\cong G/H$.

  We denote also  by $\pi$ the holomorphic extension of the
representation $\pi$ of $G$ to  the complexification $G_{\C}$ of $G$,
 whenever this
holomorphic extension exists. In fact, it can be shown that
 in the situations under interest
in this paper, this holomorphic extension exists 
\cite{neeb95,neeb96}.
 Then there exists a 
homomorphism $\chi_0 $ ($\chi$), $\chi_0: ~H\rightarrow \db{T}$,
($\chi : B\rightarrow \C^{\times}$), 
such that 
$ H = \{g\in G| e_g=\chi_0(g)e_0\}$ (respectively, 
$B = \{g\in G_{\C}| e_g=\chi(g)e_0\})$, where
$\db{T}$ denotes the torus $\db{T}:=\{ z\in\C| |z|=1\}$. 

For the homogeneous space $M=G/H$ of cosets $\{gH\}$, let $\lambda
:G\rightarrow G/H$ be the natural projection $g\mapsto gH$, and let
$o:=\lambda (\mb{1})$, where $\mb{1}$ is the unit element of $G$. 
 Choosing a section $\sigma :G/H\rightarrow
G$ such that $\sigma ( o )={\bf{1}} $, every element $g\in G$ can
be written down as $g=\tilde{g}(g)h(g)$, where
$\tilde{g}(g)\in G/H$ and $h(g)\in H$. Then we have
\begin{equation}\label{genn}
e_g=e^{i\alpha (h(g))}e_{\tilde{g}(g)},~
e^{i\alpha(h(g))}:=\chi_{0}(h)~.
\end{equation} Now we take into account that
$M$ also admits an embedding in  $G_{\C}/B$. We choose a local system
of coordinates parametrized by $z_g$ (denoted also simply $z$, where
there is no possibility of confusion) on  $G_{\C}/B$. Choosing a
section 
$G_{\C}/B
\rightarrow G_{\C}$ such that any element $g\in G_{\C}$ can be written
down
as $g=\tilde{g}_bb(g)$, where $\tilde{g}_b\in G_{\C}/B$, and
$ b(g)\in B$, we have
\begin{equation}e_g=\Lambda (g)e_{z_g},~
\Lambda (g):= \chi
(b(g))=e^{i\alpha(h(g))}(e_{z_g},e_{z_g})^{-\frac{1}{2}}. 
\end{equation}

 Let    us denote by
$\got{m}$
 the vector subspace of
the Lie algebra \g~  orthogonal to $\got{h}$, i.e. we have the vector space
decomposition $\g=\got{h}+\got{m}$. It can be shown that for CS-groups
 the vector space decomposition $\g=\got{h}+\got{m}$ is
Ad $H$-invariant. The homogeneous spaces $M\cong G/H$ with this
decomposition are called 
 {\em reductive spaces} (cf. \cite{nomizu})
and  {\em the
CS-manifolds  are  reductive spaces} (cf. \cite{last}).
So, {\it the tangent space to $M$ at $o$
can be identified with $\got{m}$}.
 % Recall 
 % that for CS-groups the
 % CS-representations are highest weight representations
 %  and the vector
% $e_0$ is a primitive element of the generalized  parabolic algebra
% $\got{b}$  (cf. \cite{neeb}). 

%Let us denote ${\mb{X}}:=d\pi(X), X\in \Ugc$, where    $\U$  denotes the
 % universal enveloping algebra.
 Let 
$\tilde{g}(g)=\exp X, \tilde{g}(g)\in G/H,~ X\in\got{m}$,
$e_{\tilde{g}(g)}=\exp({\mb{X}})e_0$.
 % Let us remember again   
 %  Theorem XV.1.1 p. 646 in \cite{neeb}.
Note  that $T_o(G/H)\cong\got{g}/\got{h}\cong \gc/\bar{\got{b}}\cong
(\got{b}+\bar{\got{b}})/\bar{\got{b}}\cong
\got{b}/h_{\C}$,
 where we
have a linear isomorphism $\alpha:\got{g}/\got{h}\cong
\gc/\bar{\got{b}}$, $\alpha(X+\got{h})=X+\bar{\got{b}}$ (cf. \cite{neeb1}). 
We can take
instead of $\got{m}\subset\g$ the subspace $\got{m}' \subset \gc$
complementary to $\bar{\bb}$, or the subspace of $\got{b}$
complementary to $\got{h}_{\C}$.
If we choose a local {\em canonical} system of coordinates
$\{z_{\alpha}\}$ with respect to the basis $\{X_{\alpha}\}$ in \m',
then we can introduce the vectors
\begin{equation}\label{cvect}
e_{z}=\exp(\sum_{X_{\alpha}\in\got{m}'}z_{\alpha}{\mb{X}}_{\alpha}).e_0
\in\Hi .
\end{equation}
We get
\begin{equation}\label{unu2}
e_{\sigma (z)}=\pi(\sigma (z)).e_0, ~z\in M, 
\end{equation}
 and  we  choose local coordinates in a neighborhood
${\mc{V}}_0\subset M$ of $z=0$ corresponding to $\sigma (o)=e\in G$ such that
\begin{equation}\label{doi}
e_{\sigma (z)}=N(z)e_{\bar{z}},~~
N(z)=(e_{\bar{z}},e_{\bar{z}})^{-1/2}.
\end{equation}
 Equations (\ref{cvect}), (\ref{unu2}), and (\ref{doi})
 define locally  the {\em coherent vector
 mapping}
\begin{equation}\label{cvm}
\varphi : M\rightarrow \bar{\Hi}, ~ \varphi(z)=e_{\bar{z}},  
\end{equation}
where $ \bar{\Hi}$ denotes the Hilbert space conjugate to $\Hi$.
We call the  vectors $e_{\bar{z}}\in\bar{\Hi}$ indexed by the points
 $z \in M $  {\it
Perelomov's coherent state vectors}.

\subsection{Reproducing kernel}\label{fock}

 Let us introduce the function $f'_{\psi}:G_{\C}\to \C$,
$f'_{\psi}(g):=(e_g,\psi ), g\in G, \psi\in \Hi$.
Then 
$f'_{\psi}(gb)=\chi(b)^{-1}f'_{\psi}(g), g\in G_{\C}, b\in B$,
where  $\chi$ is  the continuous homomorphism of the 
isotropy subgroup $B$ of  $G_{\C}$  
  in $\C^{\times}$.   
 {\it The coherent states realize the space
of holomorphic global sections} $
\Gamma^{\text{{hol}}}(M,L_{\chi})=H^0(M,L_{\chi})$ {\it on the}
 $G_{\C}$-{\it homogeneous line bundle
$L_{\chi}$ associated by means of the character
 $\chi$ to the  principal B-bundle}
 (cf. \cite{onofri}, \cite{neeb94}).
 The  holomorphic line bundle is $L_{\chi}:=M\times_{\chi}\C$,  
 also denoted   $L:=M\times_B\C$ (cf. \cite{bott,tirao}).

The local trivialization of the line bundle $L_{\chi}$ associates to
every $\psi\in\Hi$ a holomorphic function $f_{\psi}$   on a  open
set  in $M \hookrightarrow G_{\C}/B$. Let the notation $G_S:=G_{\C}
\setminus S$,
 where $S$ is the   set
$S:=\{g\in G_{\C}|\alpha_g=0\}$,
and $\alpha_g:= (e_g,e_0)$. $G_S$ is a dense subset of $G_{\C}$.
We  introduce  the function $f_{\psi}:G_{S}\mapsto\C$,
$f_{\psi}(g)=\frac{f'_{\psi}(g)}{\alpha_g}, \psi\in\Hi ,~g\in G_{S}$.
The
function $f_{\psi}(g)$ on $G_S$
 is actually a function of the natural projection
$\lambda (g)$, $\lambda : G\rightarrow G/H$, 
 holomorphic in  $M_S:=\lambda ( G_S)$.

Supposing that
{\it the line bundle $L_{\chi}$ is already very ample}, the symmetric
Fock space 
\fl~ is defined as the set of functions corresponding to sections such
that
$\{f\in L^{2}(M,L)\cap{\mc{O}}(M,L)
|(f,f)_{\fl}<\infty \}$
 with respect to  the scalar product
\begin{equation}\label{scf}
(f,g)_{\fl} =\int_{M}\bar{f}(z)g(z)d\nu_M(z,\bar{z}),
\end{equation}
where $d{\nu}_{M}(z,\bar{z})$ is the quasi-invariant measure on $M$
\begin{equation}\label{scf1}
d{\nu}_{M}(z,\bar{z})=\frac{\Omega_M(z,\bar{z})}{(e_{\bar{z}},e_{\bar{z}})}.
\end{equation}
Here  $\Omega_M$ is the $G$-invariant volume form
\begin{equation}
\Omega_M:=(-1)^{\binom {n}{2}}\frac {1}{n!}\;
\underbrace{\omega\wedge\ldots\wedge\omega}_{\text{$n$ times}}\ ,
\end{equation}
and the K\"ahler two-form $\omega$ on $M$ is given by
\begin{equation}\label{twoform}
\omega (z)
 = i\sum_{\alpha ,\beta\in\Delta_{\got{m}'}}
\NC_{\alpha,\beta}dz_{\alpha} \wedge
d\bar{z}_{\beta},~
\NC_{\alpha,\beta}(z)=\frac{\pa^2}{\pa z_{\alpha}\pa\bar{z}_{\beta}}\log
(e_{\bar{z}},e_{\bar{z}}). 
\end{equation}

It can be shown (cf. \cite{neeb94})  that the space of functions
 $\fl$ identified with $ L^{2,{\text{hol}}}(M,L_{\chi})$ {\em  is a closed
subspace of $L^2(M,L_{\chi})$ with continuous point evaluation} and
 eq.  (\ref{scf})
is nothing else than the {\em {Parseval   overcompletness
 identity}}  \cite{berezin}
 \begin{equation}\label{orthogk}
  (\psi_1,\psi_2)=\int_{M=G/H}(\psi_1,e_{\bar{z}})(e_{\bar{z}},\psi_2)
 d{\nu}_{M}(z,\bar{z}),~ (\psi_1,\psi_2\in \Hi ) .
 \end{equation}

It can be seen 
that {\it the  relation  (\ref{scf}) (or eq. (\ref{orthogk})) 
 on homogeneous
manifolds fits into 
 Rawnsley's global realization \cite{raw} of Berezin's coherent states on 
quantizable K\"ahler manifolds} \cite{berezin1},
modulo   Rawnsley's ``epsilon'' function \cite{raw,CGR1},
 a constant for homogeneous quantization. 
If $(M,\omega )$ is a K\"ahler manifold and $(L,h, \nabla )$ is a (quantum)
holomorphic line bundle $L$ on $M$, where 
$h$ is the hermitian metric and $\nabla$
is the connection compatible with the metric and the complex structure,
then $h(z,\overline{z})=(e_{\bar{z}},e_{\bar{z}})^{-1}$ and the
K\"ahler potential is $-\log h(z)$.

Let us  now introduce the map
\begin{equation}\label{aa}
\Phi :\Hi^{\star}\rightarrow \fl ,\Phi(\psi):=f_{\psi},
f_{\psi}(z)=\Phi(\psi )(z)=(\varphi (z),\psi)_{\Hi}=(e_{\bar{z}},\psi)_{\Hi},~
z\in{\mathcal{V}}_0,
\end{equation}
where we have identified the space $\overline{\Hi}$  complex conjugate 
 to \Hi~  with the dual
space
$\Hi^{\star}$ of $\Hi$.
{\em{ Our supposition that $L_{\chi}$ is already a very ample line
bundle implies the validity of  Parseval   overcompletness
 identity }}  (\ref{orthogk}) 
(cf.  Theorem XII.5.6 p. 542 in \cite{neeb},
 Remark VIII.5 in \cite{neeb94}, and  Theorem XII.5.14 p. 552 in \cite{neeb}).

It can be defined a function  $K$,
   $K: M\times\overline{M}\rightarrow \C$, which on  ${\mathcal{V}}_0\times
\overline{\mathcal{V}}_0$ reads
\begin{equation}\label{kernel}
K(z,\overline{w}):=K_w(z)=
(e_{\bar{z}},e_{\bar{w}})_{\Hi}.
\end{equation}

Taking into account (\ref{aa}),
 it follows (see
Proposition 1 in \cite{sbctim}) that 
 {\it if the line bundle $L$ is very ample,
then the function $K$ (\ref{kernel}) is a reproducing
kernel, the \FSB~\fl~
 is the reproducing kernel Hilbert space $\Hi_K\subset \C^M$
associated to the kernel $K$}, and  
 {\it the evaluation  map $\Phi$ defined in
 eqs. (\ref{aa}) 
extends to an  isometry   }
\begin{equation}\label{anti}
(\psi_1,\psi_2)_{\Hi^{\star}}=(\Phi (\psi_1),\Phi
(\psi_2))_{\fl}=(f_{\psi_{1}},f_{\psi_{2}})_{\fl}=
\int_M\overline{f}_{\psi_1} (z)f_{\psi_2}(z)d\nu_M(z) .
\end{equation}

\section{Representations of coherent state Lie algebras  by
 differential operators}\label{DIFF}
\subsection{Multipliers and  coherent states}\label{mult}

Recalling the  definition  of the function 
$f_{\psi}$
 given in \S \ref{fock}, we have
\begin{equation}\label{iar}
f_{\psi}(z)=(e_{\bar{z}},\psi)=\frac{(\pi(\bar{g})e_0,\psi)}
{(\pi(\bar{g})e_0,e_0)},~z\in M, \psi\in\Hi .
\end{equation}
We get 
\begin{equation}\label{iar1}
 f_{\pi(\overline{g'}).\psi}(z)= \mu
(g',z)f_{\psi}(\overline{g'}^{-1}.z),
\end{equation}
where 
\begin{equation}\label{iar2}
 \mu (g',z)=
\frac{(\pi(\overline{g'}^{-1}\overline{g})e_0,e_0)}{(\pi(\overline{g})e_0,e_0)}
=\frac{\Lambda (g'^{-1}g)}{\Lambda (g)},
\end{equation}
or 
\begin{equation}\label{iar3}
\mu (g',z)=\Lambda(\bar{g'})(e_{\bar{z}},e_{\bar{z'}})=e^{i\alpha
  (\bar{g'})}
\frac{(e_{\bar{z}},e_{\bar{z'}})}{(e_{\bar{z'}},e_{\bar{z'}})^{1/2}}.
\end{equation}

The following assertion is easily checked up using successively
eq. (\ref{iar2}):
\begin{Remark}Let us consider the relation (\ref{iar}). Then we have
(\ref{iar1}), where  $\mu$ can be written down as in equations
(\ref{iar2}), 
(\ref{iar3}). 
We have the relation
 $\mu (g,z) =J(g^{-1},z)^{-1}$, i.e. the multiplier $\mu$
 is the cocycle in the  unitary representation
$(\pi_K,\Hi_K)$ attached to the
positive definite holomorphic kernel $K$ defined by equation (\ref{kernel}),
 \begin{equation}\label{num}
(\pi_K(g).f)(x):=J(g^{-1},x)^{-1}.f(g^{-1}.x),
\end{equation}
and the cocycle verifies the relation
%\begin{equation}\label{prod}
$J(g_1g_2,z)=J(g_1,g_2z)J(g_2,z).$
%\end{equation}
\end{Remark}
  Note that
{\em   the prescription (\ref{num})
 defines
a continuous action of $G$ on ${\mbox{\rm{Hol}}}(M,\C )$ with respect to
the compact open topology on the space  ${\mbox{\rm{Hol}}}(M,\C )$.
 If $K:M\times \bar{M}\rightarrow \C$ is a continuous positive definite kernel
holomorphic in the first argument satisfying
$K(g.x,\overline{g.y})=J(g,x)K(x,\bar{y})J(g,y)^*$,
$g\in G$, $x,y\in M$, then the action of $G$ leaves the reproducing kernel
Hilbert space $\Hi_K\subseteq {\mbox{\rm{Hol}}}(M,\C )$ invariant and defines
a continuous unitary representation $(\pi_K,\Hi_K)$ on this space}
(cf. Prop. IV.1.9 p. 104 in  Ref. \cite{neeb}). %$\gata$

%Note also that in the notation of Perelomov \cite{per} $e^{i\beta
%(g,z)}\equiv \mu (g,z)\equiv J(g^{-1},z)$.

\subsection{Hermitian representations and differential operators}\label{sheaf}

Let us consider again the triplet $(G, \pi, \Hi )$.
 Let \g~ be the Lie algebra of $G$ and let us denote
 by $\mathcal{S}:=\Ugc$ the semigroup associated with
 the universal enveloping  algebra equipped with the antilinear involution
 extending the antiautomorphism $X\mapsto X^*:=-\bar{X}$ of $\gc$.
  {\it The derived representation}
 $d\pi$ defined by eq. (\ref{derived}) is a {\it hermitian representation}
 of $\mathcal{S}$ on $\Hi^0:=\Hinf$ (cf. Neeb \cite{neeb}, p. 30). As we have
already noted, the unitarity and the continuity of the representation
$\pi$ implies that $i d\pi (X)|_{\Hinf}$ is essentially selfadjoint. 
Let us denote his image in $B_0(\Hi^0)$ with $\am := d\pi(\mathcal{S})$. 
If $\Phi : \sHi\rightarrow \fl $ is the % (Segal-Bargmann) 
isometry
(\ref{aa}), we are interested in the study of the image of 
\am~  via $\Phi$  as subset in the
algebra of holomorphic, linear differential operators,  
$ \Phi\am\Phi^{-1}:={\db{A}}_M\subset\got{D}_M$.
%The new results for semisimple Lie groups $G$ for
% ${\db{A}}_M$, $M\approx G/H$,
% are contained in our main  Theorem \ref{bigtheorem}.

The {\it  sheaf}
 $\got{D}_M$ (or simply \D ) {\it of holomorphic, finite order, linear
differential operators on} $M$ is a
 subalgebra of homomorphisms ${\mc Hom}_{\C}({\mc O}_M,{\mc O}_M)$ 
 generated
 by the sheaf ${\mc O}_M$ of germs of holomorphic functions of $M$ and the
 vector fields. 
 We consider also {\it  the subalgebra} \AM~ of ${\db{A}}_M$~
 {\it of differential operators with
 holomorphic polynomial coefficients}.
Let $U:=\mathcal{V}_0$ in $M$, endowed with the 
coordinates
$(z_1,z_2,\cdots ,z_n)$. We set $\pa_i :=\frac{\pa}{\pa z_i}$ and
$\pa^{\alpha}:=\pa^{\alpha_1}_1 
\pa^{\alpha_2}_2\cdots \pa^{\alpha_n}_n$, $\alpha :=(\alpha_1,
\alpha_2 ,\cdots ,\alpha_n)\in\N^n$. The sections of \DM~ on $U$ are
$A:f\mapsto \sum_{\alpha}a_{\alpha}\pa^{\alpha}f$,
$a_{\alpha}\in\Gamma (U, {\mc{O}})$, the $a_{\alpha}$-s being zero
except a finite number. 

For $k\in\N$, let us denote by $\D_k$ the subsheaf of differential
operators of degree $\le k$ and by $\D'_k$ the subsheaf of elements of 
$\D_k$ without constant terms.  $\D_0$ is identified with
$\mc{O}$ and $\D'_1$ with the sheaf of vector fields. The
filtration of \DM~ induces a filtration on  ${\A}_M$.

Summarizing, we have a correspondence between the following three objects:
\begin{equation}\label{correspond}
\g \ni X \mapsto\mb{X}\in\am\mapsto\db{X}\in\AAA\subset \DM, {\mbox{\rm
 ~ {differential~operator~on}}}~ \fl .
\end{equation}

Using  eq. (\ref{orthogk}) and the reproducing kernel properties, it is
easy to  emphasize the correspondence between the operators ${\mb{L}}\in
B_0(\Hi_0)$ and their images ${\db{L}}=\Phi {\mb{L}}\Phi^{-1}$ in $\AAA$
 defined on $M_S$.

\begin{Remark}\label{rem11}
 Let us consider $\phi ,\psi\in \Hi$, and
${\mb{L}}\in B_0(\Hi_0)$ related by
\begin{equation}\label{phipsi}
\phi ={\mb{L}}\psi .
\end{equation} Then their images $f_{\phi}, f_{\psi}
\in \fl $  are related by
 \begin{equation}\label{diffaction}
f_{\phi}(z)= {\db L}(z)f_{\psi}(z),
\end{equation} where the  operator ${\db L}= \Phi {\mb{L}}\Phi^{-1}\in
\db{A}_M $ is
 determined by its symbol  $K^{\mb{L}}$, expressed locally as 
 \begin{equation}\label{nota}
K^{\mb{L}}(z,\bar{w}):=(e_{\bar{z}},{\mb{L}} e_{\bar{w}})=
{\db L}(z)(\e_{\bar{z}},\e_{\bar{w}}) .
\end{equation}
\end{Remark}

Now we can see that
 \begin{Proposition}\label{p5}
If $\Phi$ is the isometry
(\ref{aa}), then 
$\Phi d\pi(\g_{\C})\Phi^{-1}\subseteq \D_1$.
\end{Proposition}
\begin{proof} Let us consider an element in $\g_{\C}$ and his image in \DM ,
via the correspondence (\ref{correspond}), i.e.:
$$\g_{\C}\ni X \mapsto{\db{X}}\in \DM;~~
 {\db{X}}_z(f_{\psi}(z))= 
\db{X}_z(e_{\bar{z}},\psi)= (e_{\bar{z}},\mb{X}\psi).$$
%where 
%$$\mb{X}:=d\pi (X):=\frac{d}{dt }|_{t=0}\pi(\exp (tX)).$$
The action $G\times M\rightarrow M$  is a holomorphic one and we
have successively:
\begin{eqnarray*}
\db{X}_z(f_{\psi}(z)) & = & 
(e_{\bar{z}}, d\pi (X)\psi )
=\frac{d}{dt}|_{t=0}(e_{\bar{z}}, \pi (\exp (tX))\psi) \\
 & =  & \left.\frac{d}{dt}\right|_{t=0} f_{\pi (\exp (tX))\psi }(z)= 
\left.\frac{d}{dt}\right|_{t=0}\mu ( \exp (tX),z)f_{\psi}(\exp (-t X).z)\\
& = & \left.\frac{d}{dt}\right|_{t=0}\mu ( \exp (tX),z))f_{\psi}(z)+\mu (0 ,z)
 \left.\frac{d}{dt}\right|_{t=0}f_{\psi}(\exp (-t X).z).
\end{eqnarray*}
We have finally
\begin{equation}\label{sss}
\db{X}_{z}(f_{\psi}(z))=\left(P_{X}(z)+
\sum Q^i_X(z)\frac{\pa}{\pa  z_i}\right)f_{\psi}(z),
\end{equation}
where $$P_X(z):=\left.\frac{d}{dt}\right|_{t=0}\mu (\exp (tX),z),~
 Q^i_X(z):=\left.\frac{d}{dt}\right|_{t=o}(\exp(-tX).z)_i.~~
$$
\end{proof}
Now we formulate the following assertion:
\begin{Remark}\label{main}
If $(G,\pi)$ is a CS-representation, then \AAA~ is a subalgebra of
holomorphic differential operators with polynomial coefficients, i.e. 
 $\AAA \subset \AM\subset \DM$. 

More exactly, for    $X\in\g$ and
    $\mb{X}:= d\pi(X)\in\am $,  let us consider
his image  $\db{X}\in \AAA$ as in relation (\ref{correspond}),
acting  on the space of functions \fl . Then, for 
CS-representations, we have that $\db{X}\in \A_1=\A_0\oplus \A_1'$.
 
Explicitly, if  $\lambda\in\Delta $ is a root and
 $X_{\lambda}$ is in a base of the Lie algebra $\g_{\C}$ of $G_{\C}$, then
 his image  ${\db X}_{\lambda}\in\DM$ acts as a first order
 differential operator on the
\FSB \fl 
 \begin{equation}\label{generic}
{\db X}_{\lambda}= P_{\lambda}+\sum_{\beta \in\Delta_{\m '}} Q_{\lambda
 ,\beta} \pa_{\beta},~ \lambda \in \Delta,
 \end{equation}
where $ P_{\lambda}$ and $ Q_{\lambda ,\beta}$  are 
 polynomials in
 $z$ and $\m '$ is the subset of $\g_{\C}$ which appears in the
definition  (\ref{cvect}) of the coherent  vectors.
\end{Remark}
 
 Actually, we don't have a proof of this 
assertion for the general case of CS-groups.
  For the compact case, there exists the calculation of
Dobaczewski \cite{dob}, which in fact can be extended
also to real semisimple Lie algebras. For compact hermitian symmetric
spaces it was shown \cite{sbcag} that the degrees of the polynomials
 $P$ and $Q$-s are   
 $\le 2$ and similarly for the non-compact hermitian symmetric  case
\cite{sbl}.
Neeb \cite{neeb} gives a proof of this Remark
  for 
CS-representations for  the (unimodular) Harish-Chandra type
groups. Let us also remember that:
{\it If $G$ is an admissible Lie group such that the universal
complexification $G\rightarrow G_{\C}$ is injective and $ G_{\C}$ is
simply connected, then $G$ is of Harish-Chandra type}
(cf. Proposition V.3 in \cite{neeb94}). 
   Differentiating eq. (\ref{num}) in order to obtain the derived
 representation (\ref{derived}),  we get two 
 terms, one in $\D_0$ and the
other one in $\D_1'$, as was shown in Proposition \ref{p5}.
A proof that the two parts are in fact $\A_0$ and respectively $\A_1'$
    is contained  in
Prop. XII.2.1 p. 515 in \cite{neeb}
  for the groups of Harish-Chandra type in the particular situation
where the space $\got{p}^+$ in  Lemma VII.2.16 p. 241 in \cite{neeb}
 is abelian. We  present below 
 explicit formulas for semisimple  Lie groups and 
  also the 
simplest example where the maximum degree of $P$ and $Q$ is   3. 
 \hfill   $\gata$\\[2ex]

\subsection{Simple examples}\label{SE}

%%%%% \begin{example}
\subsubsection{ Canonical commutation relations,  
 Glauber's coherent states and   the Heisenberg-Weyl  Group}
\label{exemplu}

 The example of the HW group  is sketched here
 only to
check up  that the  formalization in previous sections leads in particular
 to the  standard
realization of the canonical commutation relations  (CCR) on \fl~  \cite{bar},
 i.e. ${\mb{a}}\mapsto \frac{\pa}{\pa z},{\mb a}^+\mapsto z$.      

The  HW group here is the group with the
3-dimensional real  Lie algebra isomorphic to the Heisenberg algebra
 $\got{h}_1\equiv
\g_{HW}=$ $<is1 +z{\mb a}^+-\bar{z}{\mb a}>_{s\in \R, z\in\C}$,
 where the bosonic creation (annihilation)
operators ${\mb a}^+$ (respectively ${\mb a}$) verify the CCR relations 
$ [{\mb a},{\mb a}^+]=1,$
and the action of the annihilation operator on the vacuum is
${\mb a}\e_0=0. $

 Let $\Hi :=L^2(\R, dx)$. Then $\fl :=\Gamma^{{\text{hol}}}(\C,\frac{i}{2
\pi}\exp(-|z|^2)dz\wedge d\bar{z})$.  The infinite-dimensional 
irreducible unitary
Schr\"odinger representations  
$\pi_{\lambda}$ of the HW group are indexed by $\lambda\in\R$, where the
infinitesimal character of the representation is $\chi_{\lambda} '(z)=
2\pi i\lambda$, $z\in\got{z}$, the center of the Lie algebra of the group,
 and we
take the standard representation ($\lambda = 1$). 
 The CS-manifold $M $  for the HW group is the quotient  $HW/\R\approx \C$. 
Let us choose the section $\sigma :M\approx \C\rightarrow HW, \sigma
(z)= (0,z)$.  The CS-vectors (\ref{unu2}) for the HW-group (Glauber's CS field
 \cite{gl})
 are given by the unitary displacement operator acting on the ground state
\begin{equation}\label{glauber}
e_{\sigma (z)}:=\exp (\bar{z}{\mb a}^+-z{\mb a})e_0=
e^{-\frac{|z|^2}{2}}e_{\bar{z}},
\end{equation}
where the Perelomov's CS-vectors  are
\begin{equation}\label{basic}
\e_z:=\exp(z{\mb a}^+)\e_0,
\end{equation}
and the constant $N$ of eq. (\ref{doi}) here  has the value
 given in (\ref{glauber}) because
\begin{equation}\label{scp}
(\e_{z'},\e_{z})=\exp(\bar{z}'z).
\end{equation}

 The coherent  vectors are eigenvectors of the annihilation operator
%\begin{equation}\label{vp}
${\mb a}\e_z=z\e_z.$
%\end{equation}

It is easy to see that
\begin{equation}\label{simbola+}
(\e_{z'},{\mb a}^+\e_z)=\bar{z}'(\e_{z'},\e_z),
\end{equation}
which is compatible with the formal equation
\begin{equation}\label{aplus}
{\mb a}^+\e_z=\frac{\pa}{\pa z}~\e_z,
\end{equation}
a formula also noted by Glauber \cite{gl}.

Equation (\ref{scp}) leads
 to the known expression of the reproducing kernel for  $M\approx \C$
\begin{equation}\label{reprk}
K(z,\bar{w}):=(e_{\bar{z}},e_{\bar{w}}):=f_{e_{\bar{w}}}(z)=\exp(z\bar{w}).
\end{equation}
$K_{\bar{w}}:z\mapsto e^{z\bar{w}}$ are contained in \fl~ and $K$ is a
positive kernel on \C .
We find
\begin{equation}
K^{{\mb a}^+}(z,\bar{w})=(e_{\bar{z}},{\mb a}^+e_{\bar{w}})=
\frac{\pa}{\pa \bar{w}}
(e_{\bar{z}},e_{\bar{w}})=z(e_{\bar{z}},e_{\bar{w}}),
\end{equation}
i.e.
\begin{equation}\label{difa+}
\Phi {\mb a}^+ \Phi^{-1}(z)=z\in \A_0.
\end{equation}

Also
\begin{equation}
K^{\mb a}(z,\bar{w})=(e_{\bar{z}},{\mb a}e_{\bar{w}})=
\bar{w}(e_{\bar{z}},e_{\bar{w}})=
\frac{\pa}{\pa z}(e_{\bar{z}},e_{\bar{w}}),
\end{equation}
i.e. 
\begin{equation}\label{difa}
\Phi \mb{a}\Phi^{-1}(z)=\frac{\pa}{\pa z}\in \A'_1.
\end{equation}

 The
operator 
   $\Phi{\mb a}\Phi^{-1}$
($\Phi{\mb a}^+\Phi^{-1}$) corresponding to ${\mb a}$
(respectively, ${\mb a}^+$)
is acting  on the pre-Hilbert space 
 $\Hi^0_K\subset\fl$ corresponding to the reproducing
kernel $K(z,\bar{w})$ (\ref{reprk}), $w$ fixed.  
 $\mb{a}$ and ${\mb a}^+$ are formal adjoint on the pre-Hilbert
space $\Hi^0_K$, $(\mb{a}v,w)=(v,{\mb a}^+w)$, $v,w\in \Hi^0_K$, 
and ${\mb a}^+$ is
$\mb{a}^{\sharp}$ in the notation of \S \ref{DIFF}.  

 Note that  the {\em principal vectors}
$e_{\bar{w}}\in\fl$ in Bargmann's terminology  (see
eq. (1.10) in \cite{bar}),
 $e_{\bar{w}}(z)=e^{\bar{w}z}$, correspond to the  coherent vectors
(\ref{basic}) parametrized with $\bar{w}$, i.e.
$e_{\bar{w}}=\exp(\bar{w}{\mb a}^+)e_0$.
 The
isometry (\ref{aa}), $\Phi :\sHi\rightarrow \fl$, which sends one
base in another one,
 $\Phi (\frac{{{\mb a}^+}^n}{\sqrt{n!}}e_0)=\frac{z^n}{\sqrt{n!}}$,   is the
 isometry (\ref{anti}), i.e.  
 $(\psi^1,\psi^2)_{\sHi}=(\Phi(\psi^1),\Phi(\psi^2))_{\fl} =
(f_{\psi^{1}},f_{\psi^{2}})_{\fl}$.

\subsubsection{$\got{sl}(2,\C )$}

 Let us now consider  the generators of $\got{sl}(2,\C )$
\begin{equation}\label{JJ}
J_+=\left(\begin{array}{cc}0 &1\\0 & 0\end{array} \right),~
J_-=\left(\begin{array}{cc}0 &0\\1 & 0\end{array} \right),~
J_0=\frac{1}{2}\left(\begin{array} {cc}1 &0 \\0 & -1\end{array}
\right) ,
\end{equation}
which verifies the commutation relations:
$$[J_0,J_{\pm}]=\pm J_{\pm}; [J_-,J_+]=-2J_0. $$
Then we can see that
\begin{Remark}\label{remark11p} Proposition \ref{p5} for
$\got{sl}(2,\C )$ is realized as 
\begin{equation}\label{JJ2}
 {\db{J}}_+=-\frac{\pa}{\pa z},~~
{\db{J}}_-=-2jz +z^2\frac{\pa}{\pa  z},~~
{\db{J}}_0= j-z\frac{\pa}{\pa z}.
\end{equation}
\end{Remark}
\begin{proof} Indeed, let 
$$g'=\left(\begin{array}{cc} a' & b'\\c' & d' \end{array}\right)
\in\got{sl}(2,\C ) .$$
Then eq. (\ref{iar1}) becomes 
$$f_{\pi (g')\psi}(z)=(a'-c'z)^{2j}f_{\psi}[(d'z-b')(a'- c'z)^{-1}].$$
For example, the calculation for $X=J_+$ corresponds to 
$$g'=g'(t)=\left(\begin{array}{cc} 1 &  t \\ 0  & 1 \end{array}\right);~
 \mu (g',z)=1; ~ f_{\pi (g') \psi}(z)=f_{\psi}(z-t),$$
and by taking the derivatives at $t=0$ we get the first relation
(\ref{JJ2}). The other ones are obtained similarly.
\end{proof}

If we use the CS-vectors $ e_z=e^{z{\mb{J}}_+}e_{j,-j}$ for the 
minimal weight, i.e.
$${\mb{J}}_+ e_{j,-j}  \not=  0;~
{\mb{J}}_- e_{j,-j}  =  0;~
{\mb{J}}_0e_{j,-j}  =  -je_{j,-j},$$
then we get formally    
\begin{equation}\label{JJ1}
{\mb{J}}_+e_z  =  \pa e_z;~{\mb{J}}_-e_z  = (2jz-z^2\pa )e_z;~
{\mb{J}}_0e_z  =  (-j+z\pa )e_z .
\end{equation}
Equations (\ref{JJ2}) and (\ref{JJ1}) differs by  an overall
``-'' sign. See also Remark \ref{remark12}.

\section{The semisimple case}\label{sSIMPLE}

\subsection{Perelomov's coherent  vectors for semisimple 
Lie groups}\label{perv}%\footnote{Lisiecki p 252, dupa Borel; 
%Harish-Chandra IV, V;
%Enright, Howe, Wallach}

All representations of compact Lie groups are CS-representations
because these representations are highest weight representations. 
Kostant and Sternberg \cite{ks} showed that for any representation of
a compact group $G$ the orbit to a projectivized highest weight vector is the 
only K\"ahler coherent state orbit.
Harish-Chandra \cite{hc} has defined highest weight representations
for non-compact semisimple (or even reductive) Lie groups. He has classified
 square integrable highest weight representations. This classification has
been fully realized by Enright, Howe and Wallach, and independently
by  Jakobsen \cite{eww}.
  Lisiecki has emphasized  (cf. \cite{lis1} and
  Theorem 6.1 in \cite{lis}) that:
 {\em a non-compact semisimple Lie group is a CS-group if and 
only if it is hermitian. If this is the case, the
$CS$-representations of $G$ are precisely the highest weight representations. 
Each of them has a unique $CS$-orbit, which is the orbit through highest
 line}. The starting point of the proof of
Lisiecki is the paper of Borel \cite{bor}, where it is proved: {\em
a noncompact semisimple Lie group $G$ admits a homogeneous
K\"ahler manifold if and only if it is of hermitian type, and such a
manifold is of the form $G/Z_{G(S)}$, where $Z_{G(S)}$ is the centralizer of
a torus $S\subset G$; moreover,    it is a holomorphic fiber bundle over the
Hermitian symmetric space $G/K$, where $K$ is a maximal compact
subgroup of $G$,
 with (compact) flag manifolds $K/Z_{G(S)}$
as fibers}.

Let us consider again the triplet $(G,\pi,\Hi )$ where $(G,\pi)$ is a
CS-representation. Then this representation can be realized as an
extreme weight representation. For linear connected 
reductive groups with $Z_K(\got{z})=\got{k}$, where
 $\got{z}$ denotes the center of the
Lie algebra $\got{k}$ of $K$, the effective representation is
furnished by the Harish-Chandra theorem 
(cf. e.g. \cite{knapp}, p. 158).  The theorem furnishes the holomorphic
discrete series for the non-compact case, and  for the compact case it
is 
equivalent with the Borel-Weil theorem  (\cite{bw}; also
cf. \cite{knapp}, p. 143).

We use standard notation referring to Lie algebras of a complex
 semisimple Lie group $G$  
 \cite{wolf}. In this case $\Delta\equiv\Delta_s$, i.e. $\Delta_r=\{\emptyset\}$,
 i.e. all roots are semisimple.

 ${\got g}$ -- complex semisimple Lie algebra

 ${\got t}\subset\got{g}$ -- Cartan subalgebra

 ${\got b}={\got t}+{\got b}^{u}$ -- Borel subalgebra

 ${\got b}^{u}=\sum_{\alpha\in\Sigma^+}{\got g}_{\alpha}$ -- the nilradical of
 ${\got b}$

 $\Sigma$ --  root system for $({\got g} ,{\got t})$
  
$\Sigma^+$ -- a positive root system
 
$\Psi $ -- a simple root system for  $\Sigma$

 $\Sigma\ni\alpha = \sum _{\mu\in\Psi}n_{\mu}(\alpha )\mu$ -- unique,
 $ n_{\mu}\in \N,
 n_{\mu} (\alpha )\ge 0~{\mbox{\rm {if~}}} \alpha \in \Sigma^+;
 n_{\mu} (\alpha )\le 0~{\mbox{\rm {if~}}}  \alpha \in \Sigma^-$ 
%\begin{flushleft}
%\hspace*{3cm} \[
%\left\{
%\begin{array}{l}
% n_{\mu}\in \N,\\
% n_{\mu} (\alpha )\ge 0~{\mbox{\rm {if~}}} \alpha \in \Sigma^+;\\
% n_{\mu} (\alpha )\le 0~{\mbox{\rm {if~}}}  \alpha \in \Sigma^-. 
%\end{array} 
% \right.
% \]
%\end{flushleft}

  $ \Psi \supset \Phi \longrightarrow
 \Phi^r =\{\alpha \in\Sigma ; n_{\mu}(\alpha )=0$
 whenever $\mu \notin \Phi\}$

 $ \Phi^u =\{\alpha \in\Sigma ; n_{\mu}(\alpha )>0$
 for some  $\mu \notin \Phi\}=\Sigma^+\setminus\{\Sigma^+\cap\Phi^r\}$

 ${\got p}_{\Phi}={\got p}^r_{\Phi}+
{\got p}^{u}_{\Phi}$ - parabolic subalgebras of ${\got
 g}$ corresponding to $\Phi\subset\Psi$

 ${\got p}^r_{\Phi}={\got t}+\sum_{\alpha\in\Phi^r}{\got g}_{\alpha}$ -- 
the reductive part of ${\got p}_{\Phi}$

  ${\got p}^{u}_{\Phi}=\sum_{\alpha\in\Phi^u}{\got g}_{\alpha}$
-- the unipotent part of ${\got p}_{\Phi}$ 

 $\Delta_0=\Phi^r; ~\Delta_-=-\Phi^u; ~\Delta_+=\Phi^u$

 $B=\{g\in G; \Ad (g){\got b}={\got b}\}$ -- Borel subgroup
 (maximal solvable)

 $P=\{g\in G; \Ad(g){\got p}={\got p}\}$ --  parabolic subgroup (contains  a
 Borel subgroup).

In the notation of Definitions
 VII.2.4 p. 234, VII.2.6 p. 236 and VII.2.22, p. 244 in \cite{neeb}
 we have $\Phi^u \equiv\Delta^+_p$
 and
$\Phi^r\equiv\Delta_k$.

We  also need  the commutation relations in the Cartan-Weyl basis \cite{helg}
\begin{equation}\label{cw}
 \left\{
 \begin{array}{ccll}
 \left[ H_i,H_j \right]   & = &   0, & i=1,\dots ,r, H_i\in\got{t},\\
  \left[ H_i,E_{\alpha}\right]  & = &\alpha_i E_{\alpha}, &  \alpha_i =
\alpha (H_i) ,\\
  \left[ E_{\alpha},E_{\beta}\right]  & =& 
 n_{\alpha ,\beta}E_{\alpha +\beta}, & \alpha +\beta \in\Delta\setminus
 \{0\}  ,\\
\left[ E_{\alpha},E_{\beta}\right]  & =& 
 0, & \alpha +\beta \notin\Delta\cup
 \{0\},\\
 \left[ E_{\alpha} ,E_{-\alpha}\right]  & =&
 H_{\alpha}=\sum\alpha_iH_i. & 
\end{array} 
\right.
 \end{equation}
As a consequence, we have also  the commutation relations:
  \begin{equation}\label{cw1}
 \left\{
 \begin{array}{l}
 \left[E_{-\gamma},E_{\gamma}\right]=-\gamma H,~
\gamma H : =  (\gamma , H)  = \sum ^r_{j=1}\gamma_jH_j; \\
 \left[H,E_{\alpha}\right]=\alpha (H) E_{\alpha}.
\end{array}
 \right.
\end{equation}

If the extreme weight $j$ (here minimal) of the representation has the 
components
 $j = (j_1,\cdots ,j_r)$, 
  where $ r$ is the  rank 
 of the Cartan algebra,  then
 \begin{equation}\label{acth}
 \left\{
 \begin{array}{l}
 \mb{H}_k \e_{j}  =  j_k \e_{j}, k=1,\dots , r ;\\
 \mb{E}_{\alpha}\e_{j}  =  0, \alpha \in \Delta_-\cup  \Delta_0.
\end{array}
 \right.
\end{equation}

Now we take into account that $(\pi , \Hi )$ is a {\bf unitary
representation} of the group $G$ on the Hilbert space $\Hi$. Recall that
$id\pi (X)|_{\Hinf}$ is {\bf essentially selfadjoint}. If 
$\{ H_k, E_{\alpha} \}$
is the Cartan-Weyl base (\ref{cw}) of complex Lie algebra $\got{g}$,
 a base of the compact real form of $\got{g}$ is $iH_k, ~
i(E_{\alpha}+E_{-\alpha}), ~E_{\alpha}-E_{-\alpha}$,
$\alpha\in\Delta_+$. The essentially  selfadjointness condition
implies that $H^*_k=H_k$ and ${\bf{E}}^*_{\alpha}={\bf{E}}_{-\alpha}$,
$\alpha\in\Delta_+$. A base of a real
(noncompact) form  of $\got{g}$
is $iH_k$, 
$E_{\alpha}+E_{-\alpha}$, $-i(E_{\alpha}-E_{-\alpha})$. If we denote
$\bf{K}_{\alpha}:= i{\bf{E}}_{\alpha}$, then we have also
$\bf{K}^*_{\alpha}=\bf{K}_{-\alpha}$. So, it is convenient to introduce the
notation
$$\bf{F}_{\alpha}:=
\left\{\begin{array}{l}{\bf{E}}_{\alpha}, ~{\text{for the compact
case}}
\\ {\bf{K}}_{\alpha}, ~{\text{for the noncompact case}}\end{array}\right. .$$
For any element $X\in\got{g}$ the corresponding ${\bf{X}}\in\am$ is a
linear combination $${\bf{X}}=\sum c'_ki{\bf{H}}_k+\sum_{\alpha\in\Delta_+}
b'_{\alpha}{\bf{F}}_{\alpha}-\bar{b}'_{\alpha}{\bf{F}}_{-\alpha},$$ 
 where $c'_k\in\R,
~b'_{\alpha}\in\C .$ So, for $G\ni g=e^X, ~X\in\got{g}$ we have
the following realization
of equation (\ref{genn})
 $$e_{g,j}=\exp ({\bf{X}})e_j=e^{i\alpha (g)}e_{b,j}, ~
e_{b,j}:=\exp (\sum_{\alpha\in\Delta_+}
b_{\alpha}{\bf{F}}_{\alpha}-\bar{b}_{\alpha}{\bf{F}}_{\alpha}) e_j.$$
In accord to (\ref{cvect}), the Perelomov's CS-vectors are
\begin{equation}\label{t}
e_{b,j}=N(z)e_{z,j},~ 
e_{z,j}=\exp (\sum_{\alpha\in\Delta_+}z_{\alpha}{\mb{F}}_{\alpha  })e_j,
  \end{equation}
where $z_{\alpha}$ are local coordinates for the coordinate
 neighborhood
  ${\mathcal{V}}_0\subset M$.

\subsection{Differential operators on semisimple Lie group orbits}\label{MAIN}

We start  introducing the   notation
 
$$ Z:=\sum_{\alpha\in\Delta_+}z_{\alpha}\mb{E}_{\alpha}, ~
\partial_{\alpha}(Z)=\mb{E}_{\alpha}, ~~ \partial_{\alpha}=
 \frac{\partial}{\partial z_{\alpha}},~ \alpha \in \Delta_+ .$$ 

With this notation, the Perelomov's coherent state vectors are
\begin{equation}\label{ceper}e_{z,j}=\exp Ze_j,
\end{equation}
but when not necessarily, the subindex $j$ will be omitted.

In this paragraph we use a formal method  to get the holomorphic
 differential
action (\ref{generic}) of a generator $X$ of the Lie algebra $\got{g}$
of the group $G$  on the homogeneous space
$M=G/H$.
 This method was developed in \cite{sbcag} (see
also \cite{morse}) and applied in \cite{sbl}.

Let us consider Perelomov's coherent state vectors (\ref{ceper}). 
We  associate to every generator $X\in\got{g}$ a
formal operator $D_{\bf{X}}$ on $\am$,
 where ${\bf{X}}:=d\pi (X)$.  Then we make the following 
\begin{Remark}\label{remark12}
Let us suppose that we have the relation
$${\bf{X}}.e_z=D_{\bf{X}}(z).e_z,$$
where $e_z$ is the Perelomov's state vector (\ref{cvect}), in
particular
 (\ref{ceper}), belonging to
the Hilbert space $\Hi$ of the unitary continuous representation
$\pi$, and $z$ are local coordinates on the homogeneous manifold
$M=G/H$. We suppose that $D_{\bf{X}}$ is a first order differential
operator with polynomial coefficients of the form
(\ref{generic}). Then the differential action $\db{X}$ on the
symmetric Fock space \fl 
\begin{equation}
\db{X}_z(e_{\bar{z}},e_{\bar{w}}):=(e_{\bar{z}},{\bf{X}}.e_{\bar{w}})
\end{equation}
is given by
\begin{equation}
\db{X}_z=D_{\bf{X}^+}(z).
\end{equation}
We can also write down  the relation
\begin{equation}\label{bidiff}
{\db{X}}_{\bar{w}}(e_{\bar{z}},e_{\bar{w}})=({\db{X}}^+)_{z}
(e_{\bar{z}},e_{\bar{w}}).
\end{equation}
\end{Remark}
\begin{proof} Indeed, let
$$f_{e_{\bar{w}}}(z)=(e_{\bar{z}},e_{\bar{w}})=K(z,\bar{w}).$$
Then
$$K^{\bf{X}}(z,\bar{w})=(e_{\bar{z}},{\bf{X}}.e_{\bar{w}})=
({\bf{X}}^+.e_{\bar{z}},e_{\bar{w}})=D_{{\bf{X}}^+}(z)(e_{\bar{z}},e_{\bar{w}}).
         $$
\end{proof}

If $G$ is a Lie group and $\got{g}$ is its Lie algebra, we
 shall use the formula (cf. \cite{bbk1} III, \S 6.4,
Corollary 3, p. 313 )
\begin{equation}\label{bch1}
{\mbox{\rm{Ad}}}(\exp Z)=\exp {\mbox{\rm{ad}}}_Z,~ Z\in\got{g},
\end{equation}
i.e. (cf. \cite{bbk1}, II, \S 6.5, eq. (22)):
 \begin{equation}\label{bch}
 \e^ZX\e^{-Z}=\sum_{n \ge o} \frac{1}{n!}\ad^n_ZX,~ X, Z\in\got{g},
 \end{equation}
where
 $${\mbox{\rm ad}}_YX=[Y,X], ~~ 
 {\mbox{\rm ad}}^m_YX=[Y,{\mbox{\rm ad}}^{m-1}_Y X],~~ m>1, ~~\ad^0_YX=X.   $$ 
 We also use the relation %(see Prop. 12, p. 314 in \cite{bbk1})
 \begin{equation}
 \e^Z\partial_{\alpha} (\e^{-Z})=-\left[\partial_{\alpha}(Z)+\sum_{n\ge 1}
 \frac{1}{(n+1)!}\ad^n_Z\partial_{\alpha}(Z)\right],
 \end{equation} 
 $$\partial_{\alpha}(Y)=\partial_{\alpha}Y-Y\partial_{\alpha}=
 -\ad_Y(\partial_{\alpha}).$$

We recall the definition of the  Bernoulli numbers $B_i$ \cite{abr}:
 
\begin{equation}
 \label{defc}\frac{x}{1-\e^{-x}}=1 +\frac{1}{2}x+\sum_{k\ge 1}(-1)^{k-1}
 \frac{B_k x^{2k}}
 {(2k)!}=\sum_{n\ge 0}c_nx^n,
 \end{equation}
 \begin{equation}\label{cb}
 c_0=1;~~c_1=\frac{1}{2}; ~~c_{2k+1}=0;~~c_{2k}=\frac{(-1)^{k-1}}{(2k)!}B_k,
 \end{equation}
 \begin{equation}\label{bernuli}
 B_1=\frac{1}{6};~~B_2=\frac{1}{30};~~B_3=\frac{1}{42};~~
 B_4=\frac{1}{300},...
 \end{equation}

 We need:
  \begin{lemma}
 Let the relation: 
\begin{equation}\label{csum}
{\frac{1}{n!}=\sum_{k=0}^{n}c_k\frac{1}{(n-k+1)!}} .
 \end{equation}
 Then  the constants $c_k$ of eq. (\ref{csum}) verifies the
 definition (\ref{defc}).
\end{lemma}
\begin{proof}
  We have successively:
 \begin{eqnarray*}
\sum_{n\ge 0}\frac{1}{n!}x^n & = & \sum_{k,n; n\ge k}c_k\frac{x^n}{(n-k+1)!}\\
  & = & \sum_{k,m\ge 0}c_k\frac{1}{(m+1)!}x^{m+k}\\
 & = & \sum_{k\ge 0}c_kx^k\sum_{m\ge 0}\frac{x^m}{(m+1)!}.
 \end{eqnarray*}
We obtain
 $$ x\e^x  =  \sum_{k\ge 0}c_kx^k [\e ^x -1]~.{\mbox{~~~~}}$$
\end{proof}  
 We need also another  formula similar to (\ref{csum}).
\begin{lemma}\label{constd}
Let the constants $d_k$ be defined by the relation: 
 
 \begin{equation}\label{dd}
 {\frac{1}{(n+2)!}=\sum_{k=0}^{n}d_k\frac{1}{(n-k+1)!}}~.
 \end{equation}

 Then the constants $c$ and $d$ are related by 
 
 \begin{equation}\label{cdd}
 {d_k= (-1)^kc_{k+1}}~.
 \end{equation}
  
\end{lemma} 
\begin{proof}
\begin{eqnarray*}\sum_{n\ge 0}\frac{1}{(n+2)!}{x^{n+2}} & = &
 \sum_{n\ge 0}\sum_{k=0}^{n}d_k\frac{x^{n+2}}{(n-k+1)!}=  \\
& = &   \sum_{k\ge 0}d_kx^{k+1}{\sum_{m\ge
 0}\frac{x^{m+1}}{(m+1)!}}.
  \end{eqnarray*}

 So
 $$(\e^x-1)\sum_{k\ge 0}d_kx^{k+1}= \e^x-x-1 ~.$$

 \begin{eqnarray*}
 \sum_{k\ge 0}d_kx^k & = &\frac{\e^x-1-x}{x(e^x-1)}\\
 & = &
 \frac{1}{x}-\frac{1}{\e^x-1} \\
 & =  & -\frac{1}{x}\sum_{n\ge 1}c_n(-x)^n\\
 & = & \sum_{k\ge 0}c_{k+1}(-1)^kx^k .
\end{eqnarray*} 
  Eq. (\ref{defc}) was used. Equation (\ref{cdd}) is proved.
  ${\mbox{~~~}}$%\gata}$
\end{proof}

 Now we formulate the main result of the present paper:
 
\begin{Theorem}\label{bigtheorem}
Let $G$ be a semisimple Lie group admitting a CS-representation $\pi$.
If  $X_{\lambda}\in\got{g}$ is  a generator
of the group $G$, then the
 corresponding   holomorphic  first-order  differential operator
${\db X}_{\lambda}$ 
associated to the derived representation $d\pi$, ${\db
X}_{\lambda}\in\D_1=\D_0\oplus 
\D'_1$, has polynomial coefficients,   ${\db X}_{\lambda}\in {\A}_1$.
 More exactly,
  \begin{equation}\label{generic1}
{\db X}_{\lambda}= P_{\lambda}+\sum_{\beta \in\Delta_+} Q_{\lambda
 ,\beta} \pa_{\beta}, \lambda \in \Delta,
 \end{equation}
where $ P_{\lambda}$ and $ Q_{\lambda ,\beta}$  are polynomials in
 $z$ on the $G$-homogeneous CS-manifold $M$. 

Explicitly, the differential operators  $\{{\db  E}_{\alpha}, {\db H}_i\}$
corresponding to the Cartan-Weyl base $\{E_{\alpha},H_i\}$  (\ref{cw}) are
as follows:

 a) For  $\alpha\in \Delta_+$,

 \begin{equation}\label{plus}
 {\db  E}_{\alpha}=\sum_{k\ge 0}^{\nu}c_k\sum_{\beta\in\Delta_+}
 p_{k\alpha\beta}(z)\partial_{\alpha +\beta},
 \end{equation}
where  the coefficients  $c_k$, related to the Bernoulli numbers by eq.
 (\ref{cb}),
 are given  by eq. (\ref{csum}). The polynomials  $ p_{k\alpha\beta},
 k\in\N, \alpha\in\Delta_+$
  are given by the
 equation: 
 \begin{equation}\label{polp}
 p_{k\alpha\beta}(z)=
 \sum_{\stackrel{\alpha_1,\cdots ,\alpha_k}{\alpha_1+\cdots +\alpha_k =
 \beta}}n_{\alpha_1\cdots \alpha_k\alpha} z_{\alpha_1}\cdots z_{\alpha _k},
~k\ge 1, 
  \end{equation}
where
\begin{equation}\label{ndef}
 n_{\alpha_1\cdots \alpha_k\alpha}=
 n_{\alpha_1,\alpha}n_{\alpha_2,\alpha +\alpha_1}\cdots n_{\alpha_k,
 \alpha +\alpha_1+\cdots +\alpha_{k-1}},~(k\ge 1,\alpha_{0}=0), 
 \end{equation}
 and $n_{\alpha\beta}, \alpha ,\beta \in \Delta_+$ are the
 structure constants of eq. (\ref{cw}), and for $k=0$ the sum
(\ref{plus}) is just $\partial_{\alpha}$.  
 
 The expression (\ref{plus}) can be put also into a form in which the Bernoulli
 numbers are explicit:

\begin{equation}\label{plus1}
  {\db E}_{\alpha}=\partial_{\alpha}+
 \frac{1}{2}\sum_{\beta\in\Delta_+}z_{\beta}n_{\beta ,\alpha}
 \partial_{\alpha +\beta}+\sum_{k\ge 1}^{\nu}\frac{(-1)^{k-1}}{{(2k)!}}B_k
 \sum_{\beta\in \Delta_+}p_{k\alpha\beta}\partial_{\alpha +\beta}.
 \end{equation}

The degree of the polynomial $p$ has the property: $
 \mbox{\rm{degree~}} p_{k\alpha\beta}\le \nu;~
 p_{k\alpha\beta} ~$ as a function of $z$  contains only
 even powers. The table below contains the values of $\nu$.

\begin{center}
 Degree $\nu$ for simple Lie algebras
 \end{center}
 $$\boxed{\begin{array}{llll} 
  A_l:\nu = l-1 & l\ge 1 & E_6: \nu = 10 & G_2: \nu =4\\
 B_l: \nu =2l-2 & l\ge 2  & E_7: \nu = 16 & \\
 C_l: \nu =2l-2 & l\ge 2  & E_8: \nu = 28    &\\
 D_l:\nu = 2l-4 & l\ge 3 & F_4: \nu = 10 & \\
\end{array}}$$ 

b) The differential action of the
  generators of the Cartan algebra is:

\begin{equation}\label{cartan}
 {\db H}= j +\sum_{\beta\in\Delta_+}\beta z_{\beta}\pa _{\beta} .
 \end{equation}

c) If $(\alpha ,j)=0$, then

\begin{equation}\label{cartan1}
 {\db E}_{\alpha}=-\sum_{\beta\in\Delta_+}n_{\beta ,-\alpha}
 z_{\beta -\alpha} \partial_{\beta} .
 \end{equation}

d) If $\gamma \in\Delta_-$ is  a simple root, then

 \begin{equation}\label{minus}
 {\db E}_{\gamma}= j \gamma z_{-\gamma} +\sum_{k\ge 0}^{\nu}d_k
 \sum_{\delta , \beta\in\Delta_+} q_{\gamma \delta}(z)
 p_{k\delta \beta}(z)\partial_{\beta +\delta},
 \end{equation}
where the coefficients $d$ are expressed through the coefficients $c$ by
 eq. (\ref{cdd}).

 The expression of the polynomials $q_{\gamma \delta},
 \gamma\in\Delta_- ,{\delta\in\Delta_+}$ is
 \begin{equation}\label{p1}
q_{\gamma \delta} = -\gamma z_{-\gamma}\delta
z_{\delta}  + \sum _{\mu\in\Delta_+}z_{\delta - \mu - \gamma} n_{\delta -\mu -
 \gamma , \gamma} z_{\mu}n_{\mu , \delta -\mu}.
 \end{equation}

 In the case of a Hermitian symmetric space
 eq. (\ref{plus}) becomes just:
\begin{equation}
{\db E_{\alpha}}=\partial_{\alpha},
 \end{equation}
 while eq. (\ref{minus}) becomes
 \begin{eqnarray*}
 -{\db E}^-_{\alpha}& = &{\db K}^-_{\alpha}=(\alpha , j)z_{\alpha}+\frac{1}{2}
 z_{\alpha}\sum_{\beta\in\Delta_n^+}(\alpha ,\beta )z_{\beta}\partial_{\beta}
  \\
& & -\frac{1}{2}z_{\alpha}\sum_{\gamma -\alpha \in\Delta_k\setminus\{ 0\}}
 n_{\gamma ,-\alpha}
 n_{\gamma -\alpha ,-\beta}
 z_{\beta + \alpha -\gamma}\partial_{\beta}
 \end{eqnarray*}

\end{Theorem}

\begin{proof}

 a) Let  $\alpha\in \Delta_+$.

We apply  the  formula (\ref{bch}):
 
\begin{eqnarray*}
\e^ZE_{\alpha}\e^{-Z} & = & \sum_{n \ge o} \frac{1}{n!}\ad^n_ZE_{\alpha}\\
 & = & 
 \sum_{k, n\ge 0}
 \frac{c_k}{(n-k+1)!}\ad^n_ZE_{\alpha} \\
& = &  \sum_{k\ge 0}{c_k}\sum_{m\ge 0}\frac{1}{(m+1)!}
 \ad^ m_Z (\ad^k_Z E_{\alpha}).
 \end{eqnarray*}

 But 
  $$\ad^k_ZE_{\alpha}=
 \sum_{\alpha_ 1,\cdots ,\alpha_k}z_{\alpha_1}\cdots z_{\alpha _k}
 n_{\alpha_1\cdots \alpha_k\alpha}E_{\alpha+\alpha_1+\cdots +\alpha_k}~,$$
 and the expression (\ref{plus})
 is obtained by successive application of  the third commutation
 relation (\ref{cw}). The sum $\alpha +\alpha_{1}+...
+\alpha_{k}$ goes until a $k=\nu$ corresponding to the largest root
(cf.
\cite{bbk}, Chapter VI, Tables pp. 250-273). So the expression
(\ref{polp}) follows.

The relation 
$$\ad^k_ZE_{\alpha}=\sum_{\beta}p_{k\alpha\beta}(z)E_{\beta +\alpha}$$
leads to
 $$\e^ZE_{\alpha}\e^{-Z}=-\sum_{k\ge 0}^{\nu}
c_k\sum_{\beta}p_{k\alpha\beta}(z)\e^Z
 \partial_{\alpha + \beta}(\e^{-Z}).$$

The relations (\ref{plus}), (\ref{plus1}) are proved.

 For example, for the $A$-series \cite{bbk}:
 $$A_l:~ \circ\li\circ\li\cdots\li\circ\li\circ$$
$${\mbox{~~~~~~~~~~~~~~~}}\alpha_1 \di \alpha_2 \di\di\hspace*{4mm} 
 \alpha_{l-1} \hspace*{7mm}
 \alpha_l$$
 
The  maximal root is: $\alpha_1+\cdots +\alpha_l$. This implies the 
  degree $\nu$ for simple Lie algebra $A_l$. Similarly for the other
cases.

b), c) The differential actions corresponding to the
 generators of the Cartan algebra
 (eq. (\ref{cartan}) and  eq. (\ref{cartan1})) were
 calculated in \cite{sbcag,sbl} using the  formula (\ref{bch}) and the
 commutation
 relations (\ref{cw}).

 d) Let $\gamma \in\Delta_-$ be simple root.
 Then $[E_{\alpha},E_{\gamma}]=n_{\alpha \gamma}E_{\alpha + \gamma}.$
It is observed that 
 $ \alpha\in\Delta_+$,  $\alpha + \gamma \in \Delta$ implies
  $\alpha +\gamma \in \Delta_-$.

Indeed, if:
$  \alpha + \gamma \in
 \Delta_-$ then $ \gamma  = -\alpha +
 \delta, \delta \in\Delta_-, -\alpha \in \Delta_-, $
 i.e. $ \gamma~$
 is not simple. But this is  not true!

 So: 

 $$\alpha + \gamma  \left\{\begin{array}{l}
\in  \Delta_0 ~,{\mbox{\rm{or}}}\\
\in  \Delta_+ ~,{\mbox{\rm{or}}}\\
 =0.
 \end{array}\right.$$

 Now we do some preliminary calculation:
% \begin{eqnarray*}
$$ [Z,E_{\gamma}]  = \sum_
 {\stackrel{\alpha\in\Delta_+}{\alpha+\gamma\in\Delta_+}}z_{\alpha}n_{\alpha
 \gamma} E_{\alpha +
 \gamma}+ 
   \sum_{\stackrel{\alpha\in\Delta_+}{\alpha+\gamma\in\Delta_0}}
 z_{\alpha}n_{\alpha \gamma}E_{\alpha + \gamma}-  z_{\gamma}\gamma H.$$
% \end{eqnarray*}

 Next 
$$ [Z,H]=\sum_{\alpha\in\Delta_+}z_{\alpha}[E_{\alpha},H]=-\sum_{\alpha \in
 \Delta_+}\alpha z_{\alpha}E_{\alpha}.$$
 Also:
 $$[Z,E_{\beta}]=\sum_{\stackrel{\mu\in\Delta_+}{\mu + \beta\in\Delta_+}}
 z_{\mu} n_{\mu \beta}E_{\mu +\beta}, ~\beta\in \Delta_0, \beta = \alpha +
 \gamma .$$
    We have used the relations (\ref{cw}), (\ref{cw1}).

 We apply again the formula (\ref{bch}):
\begin{eqnarray}\label{ss1}
\nonumber
 \e^ZE_{\gamma}\e ^{-Z} & = & \sum_{n\ge 0}\frac{1}{n!}\ad ^n_{Z} E_{\gamma}
 =\\
\nonumber  & = & E_{\gamma} + \sum_{m\ge 0}
 \frac{1}{(m+1)!}ad ^m_Z [Z,E_{\gamma}]=\\
\nonumber & = & E_{\gamma} -
 \sum_{\stackrel{\alpha\in\Delta_+}{\alpha + \gamma \in
 \Delta_+}} z_{\alpha}n_{\alpha \gamma}\e^Z\partial_{\gamma +\alpha}(\e^{-Z})+
 \\
\nonumber & + & \sum_{\stackrel{\alpha\in\Delta_+}{\alpha + \gamma \in
 \Delta_0}} z_{\alpha}n_{\alpha \gamma}E_{\alpha +\gamma}-
 \gamma z_{_\gamma}H +  R~,
\end{eqnarray}
where 
 \begin{equation}
R := \sum _{m\ge 1}\frac{1}{(m+1)!}\ad^{m-1}_Z \left[ q\right ],
\end{equation}
and
 \begin{equation}
 \label{Q}
 q :=  -\gamma z_{-\gamma}
 \sum_{ \alpha\in\Delta_+}\alpha z_{\alpha}E_{\alpha} + 
 \sum_{\stackrel{\alpha , \mu , \alpha + \gamma \in \Delta_+}{\alpha + \gamma +
 \mu \in \Delta_+}}z_{\alpha} n_{\alpha \gamma}z_{\mu}n_{\mu , \alpha + \gamma}
 E_{ \mu + \alpha + \gamma }~ . 
 \end{equation}

Changing the summation variable $\alpha\rightarrow\delta$ in the first sum in 
 the expression  (\ref{Q}) and denoting $\mu + \alpha + \gamma \rightarrow
 \delta $ in the second sum of the same expression, we get finally for $q$ the
 formula

 \begin{equation}\label{p}
 q  = \sum_{\delta\in\Delta_+}q_{\gamma \delta} E_{\delta}~,
\end{equation}
 and the formula (\ref{p1}) is proved.

  We continue to calculate $R$:
 
$$R=\sum_{m\ge 0}\frac{1}{(m+2)!}\ad^m_Z\sum_{\delta\in\Delta_+}q_{\gamma
 \delta} E_{\delta}~ .$$

 Now we use eq. (\ref{dd}).  Then 
\begin{equation}
R  =  -\sum_{k\ge 0}d_k\sum_{\beta ,\delta\in\Delta_+}q_{\gamma ,\beta}p_{k
 \delta\beta}(z)\e^Z\partial_{\delta +\beta}(\e^{-Z})~. 
\end{equation}

 So we get finally:

 \begin{eqnarray*}
\e^ZE_{\gamma}\e^{-Z} & = & E_{\gamma}-\sum_{\stackrel{\alpha\in\Delta_+}{
 \alpha +\gamma \in \Delta_+}}z_{\alpha}n_{\alpha \gamma}
 \e^Z\partial_{\gamma}(\e^{-Z})-\\
 & & -\gamma z_{-\gamma}H +\sum_{\stackrel{\alpha\in\Delta_+}{
 \alpha +\gamma \in \Delta_0}}z_{\alpha}n_{\alpha \gamma}E_{\alpha +\gamma}
 +R ~,% \Downarrow
 \end{eqnarray*}
and eq. (\ref{minus}) is proved.% $  ~~~\gata$
\end{proof}

\subsection{Examples}\label{EXA}
\subsubsection{  $M=SU(3)/\T$}

In this section we follow closely \cite{saraceno}.

 The commutation relations of the generators are:
 
\begin{equation}\label{grun}
[C_{ij}, C_{kl}]=\delta_{jk}C_{il}-\delta_{il}C_{kj}, ~1\le i, j \le 3.
 \end{equation}

Let us consider the following parametrizations useful for the Gauss
decomposition and also in the definition of the coherent states for
the manifold $M$:
\begin{equation}\label{vzeta}
V_{+}(\zeta ):=\exp (\zeta_{12}C_{12}+\zeta_{13}C_{13}+\zeta_{23}C_{23}),
\end{equation}
\begin{equation}\label{vz}
V'_{+}(z ):=\exp (z_{23}C_{23})\exp (z_{12}C_{12}+z_{13}C_{13}).
\end{equation}

Let us denote by the same letter $C_{ij}$ the $n\times n$-matrix
having all elements $0$ except at the intersection of the line $i$
with the column $j$, that is $C_{ij}=(\delta_{ai}\delta_{bj})_{1\le a, 
b\le n}$. Here  $n=3$. Then:
\begin{equation}\label{matricezeta}
V_{+}(\zeta )=
\left(
\begin{array}{ccc}
 1 & \zeta_{12} & \zeta_{13}+\frac{1}{2}\zeta_{12}\zeta_{23}\\
0 & 1 & \zeta_{23}\\
0 & 0 & 1
\end{array}
\right),
\end{equation}

\begin{equation}\label{matricez}
V'_{+}(z )=
\left(
\begin{array}{ccc}
 1 & z_{12} & z_{13}\\
0 & 1 & z_{23}\\
0 & 0 & 1
\end{array}
\right).
\end{equation}

Now observing that for 
\begin{equation}\label{egal}
z_{12}=\zeta_{12}; ~z_{13}=\zeta_{13}+\frac{1}{2}\zeta_{12}\zeta_{23};~
z_{23}=\zeta_{23}~, 
\end{equation}
we get 
\begin{equation}\label{eegal}
V_{+}(\zeta )=V'_{+}(z )~.
\end{equation}
So we have two  parametrizations  of the compact non-symmetric flag manifold
 $M=SU(3)/S(U(1))\times U(1)\times U(1))$: one in  $\zeta$, given by eq. 
(\ref{matricezeta}) and the other one in $z$, given by  (\ref{matricez}), which
 are identified using the relations (\ref{egal}).

Let us consider also the vectors
\begin{equation}\label{vect}
\begin{array}{ccl}
\phi_z & = & [V'_{+}(z )]^+\phi_{w}\\
& = & \exp (\bar{z}_{12}{{C}}_{21}+\bar{z}_{13}{{C}}_{31})
\exp (\bar{z}_{23}{{C}}_{32})\phi_{w} ~.
\end{array}
\end{equation}

$\phi_w$ is chosen as maximal weight vector corresponding to the weight
$w=(w_1,w_2,$ $w_3)$
 such that  $ j_1=\omega_1-\omega_2\ge 0, j_2=
 \omega_2-\omega_3\ge 0 $ and the lowering operators are ${\bf{C}}_{ij}, i>j$, 
while ${\bf{C}}_{ii}$ corresponds  to the Cartan algebra, i.e.
\begin{equation}\label{altele}
\left\{
\begin{array}{cll}
{\mb C}_{ij}\phi_w & \not=&  0,~ i>j~;\\
{\mb C}_{ij}\phi_w & = & 0 ,~i<j~;\\
{\mb C}_{ii}\phi_w & = & w_i\phi_w~.
\end{array}
\right.
\end{equation}
The coherent  vectors corresponding to the representation
$\pi_w$ determined by eqs. (\ref{altele})  are introduced as
\begin{equation}\label{noucoer}
e_z=\pi_{w}((V'_{+}(z ))\phi_w ~.
\end{equation}
Denoting by $Z$ the matrix
 \begin{equation}\label{zmat} 
Z=
 \left(
 \begin{array}{lll}
 1 & z_{12} & z_{13} \\
 0 & 1 & z_{23} \\
 0 & 0 & 1
 \end{array}
 \right)~,
\end{equation}
the reproducing kernel which determines the scalar product
$(e_{\bar{z}},e_{\bar{z}})$
has the expression: 
\begin{eqnarray*}
 K(ZZ^+) & = &\Delta^{j_1}_1(ZZ^+)\Delta^{j_2}_2(ZZ^+);\\
 \Delta_1 & = & 1+|z_{12}|^2+|z_{13}|^2;\\
\Delta_2 & = & (1+|z_{12}|^2+|z_{13}|^2)(1+|z_{23}|^2)
-|z_{12}+z_{13} \overline{z}_{23}|^2 .
 \end{eqnarray*}

In particular, it is observed that $z_{23}= 0$ corresponds to the manifold
%\begin{equation}
%\label{particular}
%z_{23}= 0\rightarrow 
$SU(3)/S(U(2)\times U(1))=G_1(\C^3)=\C\P^2$.
%\end{equation}

In order to compare with the scalar product for  coherent states on
$M\approx\C\P^2\approx G_1(\C^3)$ we remember that in the case of the
Grassmannian we have used in \cite{sbcag,sbl}
 a weight which here corresponds to $w_1=1,
w_2=w_3=0$
and then on $\C\P^2$ the reproducing kernel is just
$$K(ZZ^+)= \Delta_1 .$$

We underline that
the calculation given below, which will proof Lemma \ref{ceurile},
is algebraic, and we do
not use the value of the reproducing kernel.

Let us introduce the simplifying  notation
\begin{equation}
\label{notatie}
E(z):=V'_{+}(z )~.
\end{equation}

We shall find the operators $\tilde{{\mb{C}}}_{ij}$ such that
\begin{equation}\label{simplu}
E{\mb{C}}_{ij}=\tc_{ij}E~.
\end{equation}
Then
\begin{equation}\label{medie}
(e_{{z}} ,{\mb{C}}_{ij}e_{{z}})=(\phi_w ,E{\mb{C}}_{ij}E^+\phi_w)=
(\phi_w ,\tc_{ij} EE^+\phi_w)=
{\wt{{\db{C}}}}_{ij}(e_{{z}},e_{{z}})~.
\end{equation}
In the coherent state representation (\ref{noucoer})
$\db{C}$ is the differential
operator associated to the operator $\pi_w (C)$. 

\begin{lemma}\label{pregatire}
The operators $\tc_{ij}$ associated to the operators ${\mb{C}}_{ij}$ as in
(\ref{simplu}) are given by the formulas:
\begin{eqnarray*}
\tc_{11} & = & {\mb{C}}_{11}-z_{12} \frac{\pa}{\pa z_{12}}-z_{13}
 \frac{\pa}{\pa z_{13}}~,\\
\tc_{12} & = & \frac{\pa}{\pa z_{12}}~,\\
\tc_{13} & = & \frac{\pa}{\pa z_{13}}~,\\
\tc_{21} & = &
{\mb{C}}_{21}+z_{12}({\mb{C}}_{11}-{\mb{C}}_{22})-z^2_{12}
 \frac{\pa}{\pa
z_{12}}-(z_{13}-z_{12}z_{23}) \frac{\pa}{\pa z_{23}}-z_{12}z_{13}
\frac{\pa}{\pa z_{13}}~, \\
\tc_{22} & = & {\mb{C}}_{22}+z_{12} \frac{\pa}{\pa z_{12}}-z_{23} 
\frac{\pa}{\pa z_{23}}~,\\
\tc_{23} & = &  \frac{\pa}{\pa z_{23}} + z_{12} \frac{\pa}{\pa
z_{13}}~,\\
\tc_{31} & = & {\mb{C}}_{31}+
z_{23}{\mb{C}}_{21}-z_{12}{\mb{C}}_{32}+z_{13}({\mb{C}}_{11}-{\mb{C}}_{33})
-z_{12}z_{23}({\mb{C}}_{22}-{\mb{C}}_{33})-z_{13}^2 \frac{\pa}{\pa z_{13}}\\
& &
-z_{23}(z_{13}-z_{12}z_{23}) \frac{\pa}{\pa z_{23}}
-z_{12}z_{13} \frac{\pa}{\pa z_{12}}~,\\
\tc_{32} & = & {\mb{C}}_{32} 
+z_{23}({\mb{C}}_{22}-{\mb{C}}_{33})-z_{23}^2 \frac{\pa}{\pa z_{23}}
+z_{13} \frac{\pa}{\pa z_{12}}~,\\
\tc_{33} & = & {\mb{C}}_{33}+ z_{23} \frac{\pa}{\pa z_{23}}+
z_{13} \frac{\pa}{\pa z_{13}}~.
\end{eqnarray*}
\end{lemma}
\begin{proof}

 First, it is observed that
$$\frac{\pa E}{\pa z_{23}}={\mb{C}}_{23}E;~
 \frac{\pa E}{\pa z_{13}}={\mb{C}}_{13}E;
~\frac{\pa E}{\pa z_{12}}=({\mb{C}}_{12}-z_{23}{\mb{C}}_{13})E~.$$
Then  formula (\ref{bch}) is applied,
taking into account the commutation relations (\ref{grun}). One
important observation is that in the relation (\ref{vz}) the
generators in the second exponential commutes and in fact this
equation is  expressed in one-parameter subgroups. 

Another useful relation is
$$\exp (z_{23}{\mb{C}}_{23})\exp (z_{12}{\mb{C}}_{12}+z_{13}{\mb{C}}_{13})\exp
(-z_{23}{\mb{C}}_{23})=
\exp (z_{12}{\mb{C}}_{12}+(-z_{12}z_{23}+z_{13}){\mb{C}}_{13})$$
{\mbox{~~~~~~~~~~~~~}}%$\gata$
\end{proof}

%The expectation values are calculated using the differential action: 
%$$\mathcal{C}_{ij}(Z,\overline{Z})=(\Psi_Z,C_{ij}\Psi_Z)=
 %{\mb C}_{ij}\log K(ZZ^+)$$

 \begin{lemma}\label{ceurile}
The differential operators ${\db{C}}_{ij}$
associated to  the generators $C_{ij}$ are given by the formulas:
  \begin{eqnarray*}
{\db C}_{11} & = & -z_{12}\pa_{12}-z_{13}\pa_{13}+ w_1~, \\
  {\db C}_{12} & = &  \pa_{12}~,\\
 {\db C}_{13} & = & \pa_{13}~,\\
 {\db C}_{21} & = & -z_{12}^2\pa_{12}-z_{12}z_{13}\pa_{13}+
 (z_{12}z_{23}-z_{13})\pa_{23}+(w_1-w_2)z_{12}~,\\ 
 {\db C}_{22} & = & z_{12}\pa_{12}-z_{23}\pa_{23}+w_2~,\\
 {\db C}_{23} & = & z_{12}\pa_{13}+\pa_{23}~,\\
 {\db C}_{31} & = &
 -z_{12}z_{13}\pa_{12}-z^2_{13}\pa_{13}+
 \underline{(z_{12}z_{23}-z_{13})z_{23}}\pa _{23}+\\
 & & (w_1-w_3)z_{13}-(w_2-w_3)z_{12}z_{23}~,\\
 {\db C}_{32} & = & z_{13}\pa_{12}-z_{23}^2\pa_{23}+(w_2-w_3)z_{23}~,\\
 {\db C}_{33} & = & z_{13}\pa_{13}+z_{23}\pa_{23} + w_3~.
  \end{eqnarray*}
\end{lemma}
\begin{proof} The operators determined in Lemma \ref{pregatire} are
used taking into account  eqs. (\ref{altele}).
 
\end{proof}
 We have underlined the apparition of a third-degree polynomial multiplying the
 partial derivative of ${\db C}_{31}$.
 Note also the relation
${\db C}_{11}+{\db C}_{22}+{\db C}_{33}=w_{1}+w_{2}+w_{3}$.

\subsubsection{ $M=Sp(3, \R)/S(U(1)\times U(1)\times U(1))$}
 
% Now we consider the case   $Sp(3, \R)/
 %  S(U(1)\times U(1)\times U(1)$, i.e.
 This is an example of a non-symmetric, non-compact
manifold. Other simple examples can be constructed taking
  quotients of the groups
  $SO^*(6)$ or $SU(2,1)$.

Firstly,  note that:\\
 $Sp(3, \R)/
 S(U(1)\times U(1)\times U(1))=
 \underbrace{Sp(3, \R )/SU(3)}_{\mbox{\rm{Siegel bull}}}
 \times SU(3)/S(U(1)\times U(1)\times U(1)) .$

 The expression of the reproducing kernel is:
 $$K(\zeta\zeta^+)=\Delta^{j_1}_1(\zeta\zeta^+)
 \Delta^{j_2}_2(\zeta\zeta^+)\Delta^{j_3}_3(\zeta\zeta^+)~, $$
where
 $$\zeta = Z W; ~~W W^+ = (1-SS^+)^{-1}; S= S^t,$$
 $W$ is a $  3\times 3$  triangular matrix and
 $$\zeta\zeta^+ = Z(1-SS^+)^{-1}Z^+ ~.$$

 The differential action of the generators % $S=(s_{ij})_{i,j=1-3}$
 is given by the formulas $(i,j=1-3)$:
 \begin{eqnarray*}
 \widetilde{\db C}_{ij}
 & = &{\db C}_{ij}+\sum_{r=1}^3s_{ir}(\frac{\pa}{\pa s_{jr}} +
 \frac{\pa}{\pa s_{rj}})~,\\
\widetilde{\db X}_{ij} & = & \frac{1}{2}(\frac{\pa}{\pa s_{ij}} +
 \frac{\pa}{\pa s_{ji}})~,\\
 \widetilde{\db Y}_{ij}
 & = &\frac{1}{2}\sum_{r,r'=1}^3s_{ir}s_{jr'}
 (\frac{\pa}{\pa s_{rr'}}
 + \frac{\pa}{\pa s_{r'r}})+ 
 \frac{1}{2}\sum_{r=1}^3(s_{ir}{\db C}_{jr}+s_{jr}{\db
 C}_{ir})~.%~~~~~~~~~
%~~~~~~\gata
\end{eqnarray*} 
\subsection*{Acknowledgments}
 I would like to thank the organizers of the
Conference {\it Operator Algebras and Mathematical Physics} in Sinaia,
Romania, June 26 - July 04, 2003,
for the opportunity to present this subject at the conference. 
 I am grateful to Karl-Hermann Neeb for several useful remarks
on  a preliminary version of
the manuscript and  to Martin Schlichenmaier for
discussions and suggestions. Preliminary results on this subject have
been presented in several places of which I would like to mention
the  {\it  XVIII Workshop on Geometric Methods
 in Physics} in 
  Bia\l owie\.{z}a,
 Poland,  
   the Technische Universit\"{a}t, Darmstadt, Germany and 
the {\it Sophus Lie seminary}
 in  Berlin, Germany,
Universit\'e Libre de Bruxelles, Belgium, Universit\'e de Lille,
the {\it International Workshop on Wavlets, Quantization and Differential
Equations:
Theory and Applications} at the University of Havana, Cuba,
  Istituto de Matematicas, Unidad Cuernavaca, Mexico,
    expressing special  thanks to
A. Odzijewicz, K.-H. Neeb, T. Friedrich, M. Cahen, G. M. Tuynman, S.
 T. Ali  and C. Villegas Blas.  The author would like to thank the
referee for constructive suggestions. 

%%%%%\today 
\end{document}